\documentclass{amsart}

\usepackage{amsmath}
\usepackage{amscd}
\usepackage{amssymb}

\newcommand{\bk}{{\bf k}}
\newcommand{\bC}{{\Bbb C}}
\newcommand{\bR}{{\Bbb R}}
\newcommand{\bZ}{{\Bbb Z}}

\newcommand{\cF}{{\mathcal F}}
\newcommand{\cH}{{\mathcal H}}
\newcommand{\cV}{{\mathcal V}}
\newcommand{\cGV}{{\mathcal G}{\mathcal V}}
\newcommand{\fh}{{\mathfrak h}}

\newtheorem{theorem}{Theorem}[section]
\newtheorem{theorem/definition}{Theorem/Definition}[section]

\newtheorem{corollary}{Corollary}[section]

\DeclareMathOperator{\codim}{codim}
\DeclareMathOperator{\rank}{rank}
\DeclareMathOperator{\sign}{sign}
\DeclareMathOperator{\Res}{Res}
\DeclareMathOperator{\Ind}{Ind}

\theoremstyle{remark}

\theoremstyle{definition}
 \newtheorem{example}{Example}[section]

\begin{document}
\title
{Delocalized equivariant cohomology of symmetric products}
\author{Jian Zhou}
\address{Department of Mathematics\\
Texas A\&M University\\
College Station, TX 77843}
\email{zhou@math.tamu.edu}
\begin{abstract}
For any closed complex manifold $X$,
we calculate the Poincar\'{e} and Hodge polynomials of
the delocalized equivariant cohomology $H^*(X^n, S_n)$ 
with a grading specified by physicists. 
As a consequence,
we recover a special case of a formula for the elliptic genera of 
symmetric products in Dijkgraaf-Moore-Verlinde-Verlinde \cite{Dij-Moo-Ver-Ver}.
For a projective surface $X$,
our results matches 
with  the corresponding formulas for
the Hilbert scheme of $X^{[n]}$.
We also give geometric construction of an action of
a Heisenberg superalgebra on $\sum_{n \geq 0} H^{*,*}(X^n, S_n)$,
imitating the constructions for equivariant K-theory by Segal \cite{Seg} and
Wang \cite{Wan}.
There is a corresponding version for $H^{-*, *}$.

\end{abstract}
\maketitle

\section{Introduction}
Given a smooth manifold $M$ 
and a finite group $G$ of diffeomorphisms,
one can study the orbifold $M/G$ by 
its de Rham cohomology $H^*(M/G)$ following Satake \cite{Sat}.
It is easy to see that 
\begin{eqnarray} \label{orbifoldcohomology}
H^*(M/G) \cong H^*(M)^G.
\end{eqnarray}
In particular,
given a manifold $X$,
the permutation group $S_n$ acts on the $n$-fold Cartesian product 
$X^n$ by permuting the factors.
The orbifold $X^{(n)} = X^n/S_n$ is called 
the {\em $n$-th symmetric product} of $X$.
From (\ref{orbifoldcohomology}),
one obtains 
\begin{eqnarray} \label{orbifoldcohomology2}
H^*(X^{(n)}) \cong (H^*(X)^{\otimes n})^{S_n} = S^n(H^*(X)),
\end{eqnarray}
where $S^n(H^*(X))$ is the $n$-th graded symmetric product of $H^*(X)$.
From this
Macdonald \cite{Mac1} obtained the following formula: 
\begin{eqnarray} \label{gfPoincare1}
\sum_{n \geq 0} P_t(X^{(n)}) q^n 
= \frac{\prod_{d \; odd} (1 + t^d q)^{b_d(X)}}
{\prod_{d \; even} (1 - t^d q)^{b_d(X)}},
\end{eqnarray} 
where $P_t(X^{(n)})$ is the Poincar\'{e} polynomial of $X^{(n)}$,
and $X^{(0)}$ is a point.
In particular,
taking $t = -1$,
one obtains the the generating functional for the Euler numbers:
\begin{eqnarray} \label{gfEuler1}
\sum_{n \geq 0} \chi(X^{(n)}) q^n = \frac{1}{(1 - q)^{\chi(X)}}.
\end{eqnarray}
An alternative proof of (\ref{gfEuler1}) was given in Zagier \cite{Zag}.
It exploits the following idea:
given any orbifold $M/G$,
for $g \in G$,
consider the Lefschetz number 
$L_g = tr\; (g_*|_{H^{even}(M)}) - tr\; (g_*|_{H^{odd}(M)})$.
By standard character theory,
\begin{eqnarray*}
\chi(M/G) = \chi(H^*(M)^G) = \frac{1}{|G|} \sum_{g \in G} L_g
= \sum_{[g] \in G_*} \frac{1}{|Z_g|} L_g,
\end{eqnarray*}
where  $|G|$ is the number of elements in $G$,
$Z_g$ is the centralizer of $g$,
and $G_*$ denotes the set of all conjugacy classes of $G$.
Applying equivariant Atiyah-Singer theorem to the de Rham complex,
one gets $L_g = \chi(M^g)$,
where $M^g = \{ x\in M: g(x) = x \}$.
Hence
\begin{eqnarray} \label{Euler1}
\chi(M/G) = \frac{1}{|G|} \sum_{g \in G} \chi(M^g) 
= \sum_{[g]} \frac{1}{|Z_g|} \chi(M^g).
\end{eqnarray}
Such an approach was used by Hirzebruch to compute the signatures of symmetric
products.
See Zagier \cite{Zag}, \S 9 for details.

Unfortunately $H^*(M/G)$ is 
not always sufficient for studies on orbifolds.
For example,
it is not suitable for the equivariant $K$-theory.
To generalize the isomorphism between the K-theory and cohomology
on compact spaces given by the Chern character, 
Baum and Connes \cite{Bau-Con} defined the delocalized equivariant cohomology 
$$H^*(M, G) = (\bigoplus_{g \in G} H^*(M^g))^G$$ 
and an equivariant Chern character $ch_G: K^*_G(M) \to H^*(M, G)$
which they showed to be an isomorphism.
It turns out that $H^*(M, G)$ is also
the correct cohomology for the string theory on
orbifolds studied by 
Dixon, Harvey, Vafa and Witten \cite{Dix-Har-Vaf-Wit, Dix-Har-Vaf-Wit2}
(see also Vafa-Witten \cite{Vaf-Wit}, \S 4.1).
One easily sees that
$$H^*(M, G) \cong \bigoplus_{[g] \in G_*} H^*(M^g)^{Z_g}.$$ 
The components $H^*(M^g)^{Z_g}$ for nontrivial conjugacy classes 
correspond to the twisted sectors.
The connection between the orbifold string theory and 
the delocalized equivariant cohomology was made 
by Atiyah and Segal \cite{Ati-Seg}.
They interpreted the following formula for orbifold Euler number 
which appeared in \cite{Dix-Har-Vaf-Wit}
as the Euler number of  $H^*(M, G)$ or equivalently $K_G^*(M)$:
\begin{eqnarray} \label{Euler2}
\chi(M, G) 
= \frac{1}{|G|} \sideset{}{'}\sum_{g, h} \chi(M^{\langle g, h \rangle})
= \sum_{[g] \in G_*} \frac{1}{|Z_g|} \sum_{h \in Z_g}  \chi(M^{\langle g, h \rangle}),
\end{eqnarray}
where $\langle g, h \rangle$ is the group generated by $g$ and $h$,
the first sum is taken over commutating pairs 
$(g, h) \in G \times G$.
From (\ref{Euler1}) and (\ref{Euler2}), 
one obtains (cf. Hirzebruch and H\"{o}fer \cite{Hir-Hof})
\begin{eqnarray} \label{orbifoldEuler2}
\chi(M, G) = \sum_{[g] \in G_*} \chi(M^g/Z_g).
\end{eqnarray}
See also Roan \cite{Roa} for 
a mathematical expositions of the orbifold Euler number.
In a more recent paper by Vafa and Witten \cite{Vaf-Wit},
the following formula corresponding to (\ref{gfEuler1}) was proved:
\begin{eqnarray}
&& \sum_{n \geq 0} \chi(X^n, S_n) q^n 
= \prod_{l \geq 1} \frac{1}{(1 - q^l)^{\chi(X)}}. \label{gfEuler2}
\end{eqnarray}
The proof in \cite{Vaf-Wit} is in the spirit of
Macdonald \cite{Mac1} mentioned above.
A proof using Lefschetz numbers 
has been given by Hirzebruch and H\"{o}fer \cite{Hir-Hof}.

Vafa and Witten noticed that $\sum q^n H^*(X^n, S_n)$ 
is the Fock space of the Heisenberg superalgebra generated by
$H^*(X)$.
Motivated by this result,
Nakajima \cite{Nak} and Grojnowski \cite{Gro} independently 
obtained the geometric construction of the representation 
on $\sum q^n H_*(X^{[n]})$,
where $X^{[n]}$ is $n$-th Hilbert scheme of points of a surface $X$.
Partly motivated by a footnote in \cite{Gro},
Segal \cite{Seg} outlined the similar constructions 
for the equivariant $K$-theory $K^*_{S_n}(X^n)$.
Wang \cite{Wan} generalized Segal's constructions to the case of 
$K^*_{G_n}(X^n)$,
where $G$ is a finite group,
$X$ is a $G$-space,
and $G_n$ the wreath product of $S_n$ with $G^n$.
This leads to constructions of vertex representations via finite 
groups (see Frenkel-Jing-Wang \cite{Fre-Jin-Wan, Fre-Jin-Wan2}).

Under the isomorphism with $K_G^*(M)$ given by the equivariant Chern character,
$H^*(M, G)$ acquires a $\bZ_2$-grading.
In the physics literature,
a rule of assigning a grading possibly by fractional numbers 
has been well-known.
For a complex orbifold $M/G$,
a number $F_g$ can be defined for each component of $M^g$, $g \in G$.
Set
$$H^{p, q}(M, G) = \bigoplus_{p, q \geq 0} H^{p-F_g, q -F_g}(M^g/Z_g).$$
See e.g. Zaslow \cite{Zas} and also \S \ref{sec:delocalized}.
In this paper,
we show that for a complex manifold $X$,
$F_g$ is just half of codimension of $(X^n)^g$ for all $g \in S_n$.
One can define Poincar\'{e} and Hodge polynomials of
$H^*(X^n, S_n)$ with the above grading.
We compute such polynomials in this paper.
For algebraic surfaces,
our formulas coincide with the results in G\"{o}ttsche \cite{Got},
G\"{o}ttsche-Soergel\cite{Got-Soe} and Cheah \cite{Che}.
We also follow the constructions in Segal \cite{Seg}
and Wang \cite{Wan} to construct the representation of the Heisenberg
superalgebra on $\sum H^*(X^n, S_n)$.
We use the Poincar\'{e} duality in our construction.

Using our formula for the Hodge polynomials,
we give a mathematical proof of a special case of a formula
for the elliptic genera of symmetric products of complex manifolds
found by Dijkgraaf, Moore, E. Verlinde and H. Verlinde by physical arguments.

The same methods work for $H^{-*, *}(X^n, S_n)$.
Also most of our discussions can be carried out for wreath products.

\section{Preliminaries on the delocalized equivariant cohomology}

\label{sec:delocalized}

Let $G$ be a finite group acting on a manifold $M$,
Baum and Connes \cite{Bau-Con} defined the {\em delocalized cohomology} 
as follows:
let $\widehat{M}$ be the disjoint union of $M^g$, $g \in G$,
then there is a natural action of $G$ on $\widehat{M}$, 
set
\begin{eqnarray}
&& H^*(M, G) = H^*(\widehat{M})^G = (\oplus_{g \in G} H^*(M^g))^G, 
	\label{dec1} \\
&& H^*_c(M, G) = H^*_c(\widehat{M})^G = (\oplus_{g \in G} H^*_c(M^g))^G.
	\label{dec2}
\end{eqnarray}
Breaking into conjugacy classes, 
it is clear that
\begin{eqnarray*}
H^*(M, G) = \oplus_{g \in G_*} H^*(M^g)^{Z_g}, & 
H^*_c(M, G) = \oplus_{g \in G_*} H^*_c(M^g)^{Z_g}.
\end{eqnarray*}
Assume that $G_1$ and $G_2$ are two finite groups,
and $\phi: G_1 \to G_2$ is a group homomorphism. 
A map $f$ between a $G_1$-space $M_1$ and a $G_2$-space $M_2$
is called {\em equivariant} (with respect to $\phi$)
if $f(g \cdot x) = \phi(g) \cdot f(x)$ for any $g \in G_1$, $x \in M_1$.
Such a map induces 
an equivariant map $\hat{f}: \widehat{M}_1 \to \widehat{M_2}$,
and hence a homomorphism 
$f^*: H^*(M_2, G_2) \to H^*(M_1, G_1)$ 
and $f^*_c: H^*_c(M_2, G_2) \to H^*_c(M_1, G_1)$.
It is then easy to see that the delocalized cohomology is functorial.

If $M_1$ is a $G_1$-space and $M_2$ a $G_2$-space,
then $M_1 \times M_2$ is naturally a $G_1 \times G_2$-space. 
Furthermore,
if $g_1 \in G_1$ and $g_2 \in G_2$,
then $((M_1 \times M_2)^{(g_1, g_2)} = M_1^{g_1} \times M_2^{g_2}$,
and $Z_{(g_1, g_2)}(G_1 \times G_2) = Z_{g_1}(G_1) \times Z_{g_2}(G_2)$.
By the ordinary K\"{u}nneth theorem,
one obtains isomorphisms
\begin{eqnarray*}
&& \kappa: H^*(M_1 \times M_2, G_1 \times G_2)
\cong H^*(M_1, G_1) \otimes H^*(M_2, G_2), \\
&& \kappa_c: H_c^*(M_1 \times M_2, G_1 \times G_2)
\cong H^*_c(M_1, G_1) \otimes H^*_c(M_2, G_2).
\end{eqnarray*}

Let $M$ be a $G$-space and $G'$ is a subgroup of $G$,
we now define some functors which we will use later.
First of all, 
the identity map on $M$ is equivariant 
with  respect to the inclusion $G' \to G$.
We denote the induced homomorphism $H^*(M, G) \to H^*(M, G')$ by $\Res_{G'}^G$, 
or simply $\Res$ when there is no confusion.
With respect to the isomorphisms
\begin{eqnarray*}
H^*(M, G) \cong \oplus_{[g] \in G_*} H^*(M^{g})^{Z_{g}(G)}, &
H^*(M, G') \cong \oplus_{[g'] \in G_*'} H^*(M^{g'})^{Z_{g'}(G')},
\end{eqnarray*}
we can give $\Res$ explicitly:
if $g \in G$ is not conjugate by elements in $G$ to any element in $G'$,
then $\Res|_{H^*(M^g)^{Z_g(G)}} = 0$;
otherwise, assume that $g$ is conjugate by elements in $G$ to 
$g_1', \cdots, g'_k \in G'$ which have mutually different conjugacy classes
in $G'$,
then $H^*(M^{g})^{Z_{g}(G)} \cong H^*(M^{g'_i})^{Z_{g'_i}(G)}$
for $i = 1, \cdots, k$
and $\Res|_{H^*(M^g)^{Z_g(G)}}$ is given by
the direct sum of the inclusions 
$H^*(M^{g'_i})^{Z_{g'_i}(G)} \hookrightarrow H^*(M^{g'_i})^{Z_{g'_i}(G')}$.
We define another homomorphism $\Ind_{G'}^G: H^*(M, G') \to H^*(M, G)$ 
as follows:
for $\alpha \in H^*(M^g)^{Z_g(G')}$, $g \in G'$,
$$\Ind_{G'}^G(\alpha) = \frac{1}{|Z_g(G)|} \sum_{h \in Z_g(G)} h^*(\alpha)
\in H^*(M^g)^{Z_g(G)} \hookrightarrow H^*(M, G).$$

We now consider a $\bZ$-graing of $H^*(M, G)$ suggested by physicists.
Let $M$ be a complex manifold and $G$ acts by biholomorphic maps. 
Denote by $N^g$ the normal bundle of a component $M^g_i$ of $M^g$ in $M$.
Then there is a natural decomposition
$N^g = \oplus_j N^g(\theta_j)$,
where $N^g(\theta_j)$ is a complex subbundle on which $g$ acts as 
$e^{\sqrt{-1} \theta_j}$
where $0 < \theta_j < 2\pi$
(these angles will be called action angles).
Zaslow \cite{Zas} suggested the shift of the bigrading
on $H^{*, *}(M_i^g)$ by $(F_g, F_g)$,
and hence the shift in grading on $H^*(M_i^g)$ by $2F_g$, 
where
\begin{eqnarray} \label{shift}
F_g = \frac{\sum_j \rank_{\bC} N^g(\theta_j) \cdot \theta_j}{2 \pi}.
\end{eqnarray}
It is interesting to notice that here we are using 
some data which appeared in equivariant index theory 
(Atiyah-Singer \cite{Ati-Sin}).
Notice that $F_g$ may not be always an integer, 
but there are conditions to ensure it is.

\begin{example} \label{example:shift}
Let $X$ be a complex manifold,
then $\bZ_n$ acts on $X^n$ by cyclic permutions of the factors. 
Denote by $\sigma$ the generator of $\bZ_n$,
then $(X^n)^{\sigma} = \Delta_n(X)$, the diagonal.
Since the cyclic matrix has eigenvalues $e^{2\pi i j/n}$,
one can easily deduce that
$$N^{\sigma} = \oplus_{j=1}^{n-1} N^{\sigma}(2\pi j/n),$$
where each $N^{\sigma}(2\pi j/n)$ is isomorphic to the complex
tangent bundle of $\Delta_n(X)$.
Hence by (\ref{shift})
$$F_{\sigma} 
= \dim_{\bC} X \cdot \sum_{k=1}^{n-1} k/n = \dim_{\bC} X \cdot (n-1)/2.$$
Note that this is exactly the half of the complex codimension of the diagonal 
in $X^n$.
For this to be an integer for all $n$,
one only needs to assume that $M$ has even complex dimension.
\end{example}

\begin{example}
Let $M/G$ be a complex orbifold of dimension $n$.
Zaslow \cite{Zas} showed that $F_{g^{-1}} = \codim M^g - F_g$ 
by the following observation:
$M^g = M^{g^{-1}}$,
and if the action angles on the normal bundle of $M^g$ 
are $\{\theta_i\}$,
then those of $g^{-1}$ are $\{2\pi - \theta_i\}$.
Now if $g$ is conjugate to $g^{-1}$,
then $F_g = \frac{1}{2} \codim M^g.$
Since every element in $S_n$ is conjugate to its inverse,
this gives an alternative calculation for Example \ref{example:shift}.
\end{example}

\begin{example}
If $M$ is a Calabi-Yau manifold and $G$ is a finite automorphism group which
preserves the holomorphic volume form,
then $F_g$ is an integer for all $g \in G$.
Indeed, $g$ acts in the fiber of the normal bundle by a matrix of determinant 
$1$,
while the determinant can be computed as $\exp (2 \pi \sqrt{-1} F_g)$.
\end{example}

\section{Preliminaries on the symmetric products of graded vector spaces}
\label{sec:series}

For a $\bZ$-graded finite dimensional vector space 
(over a field $\bk \cong \bR$ or $\bC$)
$V = \sum_{d \in \bZ} V_d$ such that $b_d(V) = \dim V_d < \infty$ for all $d$
and $V_{d} = 0$ for $d < 0$,
the Poincar\'{e} series of $V$ is by definition
$$p_t(V) = \sum_{d \geq 0} b_d(V) t^d.$$
Denote by $\cGV_{\geq}$ the set of all such graded vector spaces.
Then $\cGV_{\geq}$ admits several operations.
For any integer $m$ and any graded vector space $V$,
$V[m]$ is the graded vector space with $V[m]_d = V_{d-m}$.
For positive integer $m$ and $V \in \cGV_{\geq}$,
$V[m] \in \cGV_{\geq}$.
For $V, V' \in \cGV_{\geq}$,
let $V \oplus V'$ be the graded vector space 
with $(V \oplus V')_d = V_d \oplus V'_d$.
Then $V \oplus V' \in \cGV_{\geq}$.
Alos let $V \otimes V'$ be the graded vector space with
$(V \otimes V')_d = \oplus_{p+q= d} V_p \otimes V'_q$.
Then $V \otimes V' \in \cGV_{\geq}$.
Clearly we have
\begin{eqnarray*}
&& p_t(V[m]) = t^m p_t(V), \\
&& p_t(V \oplus V') = p_t(V) + p_t(V'), \\
&& p_t(V \otimes V') = p_t(V) \cdot p_t(V').
\end{eqnarray*}
Since $\cGV_{\geq}$ is an abelian semigroup under $\oplus$,
one can consider its Grothendieck group
$GK_{\geq}$ (cf. Atiyah \cite{Ati}, \S 2.1).
Then $p_t$ is a ring homomorphism $p_t: GK_{\geq} \to \bZ[[t]]$.
Denote by $\cV$ the space of (ungraded) finite dimensional vector spaces.
This is  also an abelian semigroup,
so one can take its Grothendieck group $K$.
This is just the $K$-theory of a point.
We regard $GK_{\geq}$ as $K[[t]]$.
There is a map $D: \cV \to \bZ$ given by taking the dimensions of the vector 
spaces.
It extends to $D: K \to \bZ$.
If we regard $GK_{\geq}$ as $K[[t]]$,
then $p_t$ is nothing but the extension of $D$ to $K[[t]]$.

Now we recall some power operation in $K$-theory (Atiyah \cite{Ati}, \S 3.1).
Given $V \in \cV$, let $S^n(V)$ be the $n$-symmetric product of $V$.
Set $S^*_t(V) = \sum_{n\geq 0} t^nS^n(V)$.
This extends to a map
$S^*_t: K \to K[[t]]$.
We also regrad this as giving a map $S^*: \cV \to GK_{\geq}$:
for $V \in \cV$,
elements of $S^n(V)$ has degree $n$.
Then $D(S^*_t(V)) = p_t(S^*(V))$.
Since
$$S^n(V \oplus V') \cong \oplus_{p+q=n} S^p(V) \otimes S^q(V'),$$
for $V, V' \in \cV$, 
we get 
$$S_t^*(V \oplus V') = S^*_t(V) S^*_t(V').$$
For a one dimensional vector space $L$,
we have
$$S_t^*(L) = 1 + tL + t^2L^{\otimes 2} + \cdots,$$
hence
$$D(S_t^*(L)) = 1 + t + t^2 + \cdots = \frac{1}{1 - t}.$$
Hence by writing $V \in \cV$ as a direct sum of one-dimensional subspaces,
one has
$$p_t(S^*(V)) = D(S^*_t(V)) = \frac{1}{(1-t)^{\dim V}}.$$
One can also consider the anti-symmetric produces:
$\Lambda_t^*: K \to K[[t]]$ and $\Lambda^*: K \to \cGV_{\geq}$
given by $\Lambda_t^*(V) = \sum_{n \geq 0} t^n\Lambda^n(V)$
and $\Lambda^*(V) = \sum_{n \geq 0} \Lambda^n(V)$ for $V \in \cV$,
where elements in $\Lambda^n(V)$ are given degree $n$.
Then $D(\Lambda^*_t(V)) = p_t(\Lambda^*(V))$.
 Since
$$\Lambda^n(V \oplus V') 
\cong \oplus_{p+q=n} \Lambda^p(V) \otimes \Lambda^q(V'),$$
for $V, V' \in \cV$, 
we get 
$$\Lambda_t^*(V \oplus V') = \Lambda^*_t(V) \Lambda^*_t(V').$$
For a one dimensional vector space $L$,
we have
$$\Lambda_t^*(L) = 1 + tL,$$
hence
$$D(\Lambda_t^*(L)) = 1 + t.$$
Therefore for $V \in \cV$,
one has
$$p_t(\Lambda^*(V)) = D(\Lambda^*_t(V)) = (1+t)^{\dim V}.$$

Similarly,
for a bi-graded vector space $W = \sum_{m, n \in \bZ} W_{mn}$
such that $h_{m,n}(W) = \dim W_{mn} < \infty$ for all $m, n \in \bZ$
and $V_{mn} = 0$ if either $m < 0$ or $n < 0$,
define the Hodge series of $W$ by
$$h_{x, y}(W) = \sum_{m, n \geq 0} h_{m,n}(W) x^my^n.$$
We denote by $2\cGV_{\geq}$ the space of all such bi-graded vector spaces.
It is routine to define $W[l, m], W \oplus W', W \otimes W' \in 2\cGV_{\geq}$
for $W, W' \in 2\cGV$, $l, m \geq 0$.
Furthermore,
we clearly have
\begin{eqnarray*}
&& h_{x, y}(W[l, m]) = x^ly^mh_{x, y}(W), \\
&& h_{x, y}(W \oplus W') = h_{x, y}(W) + h_{x,y}(W'), \\
&& h_{x, y}(W \otimes W') = h_{x,y}(W) \cdot h_{x,y}(W').
\end{eqnarray*}
Denote by $2GK_{\geq}$ the Grothendieck group of the abelian semigroup
$2cGV_{\geq}$.
Then $h_{x, y}$ extends to a homomorphism 
$h_{x, y}: 2GK \to \bZ[[x, y]].$ 
We can identify $2GK$ as $GK[[y]]$ or $K[[x, y]]$.
With respect to the former identification,
$h_{x, y}$ is identified with the extension of $P_x$.
With respect to the latter, $h_{x, y}$ is identified 
with the extension of $\dim$.
 
Given $V \in \cGV$,
let $T^*(V)$ be the tensor algebra of $V$,
$I$ the ideal of $T^*(V)$ generated by elements of the form 
$v \otimes w - (-1)^{pq} w \otimes v$, 
where $v \in V_p$, $w \in V_q$.
Set $S^*(V) = T^*(V) /I$.
There is a natural decomposition
$S^*(V) = \oplus_{n \geq 0} S^n(V)$,
where $S^n(V)$ is called the {\em $n$-th graded symmetric product of $V$}.
We regard $S^*(V)$ as a $\bZ \times \bZ$-graded vector space.
For $v_1 \in V_{|v_1|}, \cdots, v_n \in V_{|v_n|}$,
denote by $v_1 \odot v_2 \odot \cdots \odot v_m$ the image of
$v_1 \otimes \cdots\otimes v_n$ in $S^n(V)$.
Then elements of such form generate $S^n(V)$,
we give such an element a bi-degree $(|v_1| + \cdots + |v_n|, n)$.
With this bi-grading,
we get a map $S^*: GK \to 2GK$ by $S^*(V) = \sum_{n\geq 0} S^n(V)$.
Alternatively,
we get a map $S^*_y: GK \to GK[[y]]$ by
$S^*_y(V) = \sum_{n \geq 0} y^n S^n(V)$.
We clearly have
$$\sum_{n \geq 0} h_{t, q}(S^n(V)) = \sum_{n \geq 0} p_t(S^n(V)) q^n.$$
In other words,
$h_{t, q}(S^*(V)) = p_t(S^*_q(V))$.
Since
$$S^n(V \oplus V') \cong \oplus_{p+q=n} S^p(V) \otimes S^q(V'),$$
for $V, V' \in \cGV$, 
we get 
$$S_t^*(V \oplus V') = S^*_t(V) S^*_t(V').$$
For a one-dimensional graded vector space $L$ which concentrates at degree $d$,
when $d$ is even
we have
$$S_q^*(L) = 1 + qL + q^2L^{\otimes 2} + \cdots,$$
hence
$$p_t(S_q^*(L)) = 1 + qt^d + q^2t^{2d} + \cdots = \frac{1}{1 - t^dq}.$$
When $d$ is odd,
we have
$$S_q^*(L) = 1 + q L,$$
hence
$$p_t(S^*_q(L)) = 1 + t^d q.$$
Hence by writing $V \in \cGV$ as 
a direct sum of one-dimensional graded subspaces,
one obtains the following well-known formula
\begin{eqnarray*}
\sum_{n \geq 0} p_t(S^n(V)) q^n 
= \frac{\prod_{d \; odd} (1 + t^d q)^{b_d(V)}}
{\prod_{d \; even} (1 - t^d q)^{b_d(V)}}.
\end{eqnarray*}
As a corollary,
one gets
\begin{eqnarray*}
\sum_{n \geq 0} p_t(S^n(V[m])) q^n 
= \frac{\prod_{d \; odd} (1 + t^d q)^{b_{d-m}(V)}}
{\prod_{d \; even} (1 - t^d q)^{b_{d-m}(V)}}.
\end{eqnarray*}

Let $W = \oplus_{p, q \in \frac{1}{2}\bZ} W_{pq}$ be 
a $\frac{1}{2}\bZ \times \frac{1}{2}\bZ$-graded vector space
such that $W^{pq} \neq \{0\}$ only if $p + q \in \bZ$.
We regard $W$ as a $\bZ$-graded vector space:
elements in $W^{p, q}$ is given the degree $p+q$.
Then the $n$-th graded symmetric symmetric product $S^n(W)$ can be defined.
It has a natural $\frac{1}{2}\bZ \times \frac{1}{2}\bZ$-grading:
for elements $w_1 \in W_{p_1q_1}, \cdots, w_n \in W_{p_nq_n}$,
$w_1 \odot\cdots \odot w_n$ has 
degree $(p_1 + \cdots + p_n, q_1 + \cdots + q_n)$.
Assume that
$h_{p, q}(W) = \dim W_{pq} < \infty$ for all $p, q \in \frac{1}{2}\bZ$,
set
$$h_{x, y}(W) = \sum_{r, s \in \frac{1}{2}\bZ} h_{r,s}(W)x^ry^s.$$
The above discussions can be repeated to obtain the following formula:
\begin{eqnarray} 
\sum_{n \geq 0} h_{x, y}(S^n(W)) q^n 
= \frac{\prod_{s+t \; odd} (1 + x^sy^t q)^{h_{s, t}(W)}}
{\prod_{s+t \; even} (1 - x^sy^t q)^{h_{s,t}(W)}}.
\end{eqnarray}
As a corollary,
one gets for $l, m \in \frac{1}{2}\bZ$ with $l + m \in \bZ$,
\begin{eqnarray} \label{gfshifted}
\sum_{n \geq 0} h_{x, y}(S^n(W[l, m])) q^n 
= \frac{\prod_{s+t \; odd} (1 + x^sy^t q)^{h_{s-l, t-m}(V)}}
{\prod_{s+t \; even} (1 - x^sy^t q)^{b_{s-l, t-m}(V)}}.
\end{eqnarray}

\section{The orbifold Poincar\'{e} polynomials of the symmetric product}

Recall that any element of $S_n$ can be uniquely written as a products of
disjoint cycles.
Denote by $N_l(g)$ the number of $l$-cycles in $g$.
The sequence $N(g) = (N_1(g), N_2(g), \cdots, )$ is called 
the {\em cycle type} of $g$.
Each cycle type corresponds to a unique conjugacy class.
Given any element $g \in S_n$ of type $N = (N_1, N_2, \cdots)$,
there is an isomorphism 
\begin{eqnarray*}
Z_g \cong S_{N_1} \times (S_{N_2} \ltimes \bZ^{N_2}_2) \times \cdots 
\times (S_{N_n} \ltimes \bZ^{N_n}_n),
\end{eqnarray*}
where each $S_{N_l}$ is given by permutating the $l$-cycles of $g$,
and 
each $\bZ_l$ is to the cyclic group generated by an $l$-cyle in $g$.
We have
$$|Z_g| = \prod_{l=1}^n N_l!l^{N_l}.$$
Denote by $\Delta_l(X)$ the diagonal in $X^l$.
It is clear that
\begin{eqnarray} \label{fix1}
(X^n)^g = \prod_{l=1}^n \Delta_l(X)^{N_l} \cong \prod_{l=1}^n X^{N_l},
\end{eqnarray}
where each cycle of $g$ contributes a factor of $X$.
Furthermore,
the action of $Z_g$ can be explicitly described as follows:
each copy of $\bZ_l$ acts trivially,  
each copy of $S_{N_l}$ permutes the $N_l$ copies of $\Delta_l(X)$.
Hence we have
\begin{eqnarray} \label{fix2}
(X^n)^g/Z_g \cong \prod_{l =1}^n X^{N_l} /S_{N_l} = \prod_{l=1}^n X^{(N_l)}.
\end{eqnarray}
Therefore, 
from (\ref{Euler1}),
one gets
\begin{eqnarray*}
&& \sum_{n \geq 0} \chi(X^{(n)}) q^n
= \sum_{n \geq 0} q^n\sum_{\sum lN_l = n} \frac{1}{\prod_{l=1}^n N_l!l^{N_l}}
\chi(\prod_{l=1}^n X^{N_l}) \\
& = & \sum_{n \geq 0} \sum_{\sum lN_l = n}\prod_{l=1}^n  
\frac{1}{N_l!}(\frac{\chi(X)q^l}{l})^{N_l}
= \prod_{l=1}^n  \sum_{N_l \geq 0}
\frac{1}{N_l!}(\frac{\chi(X)q^l}{l})^{N_l} \\
& = & \prod_{l=1}^n \exp \left( \frac{\chi(X)q^l}{l}\right)
= \exp \left( \sum_{l=1}^n \frac{\chi(X)q^l}{l}\right) \\
& = & \exp (- \chi(X) \log (1 -q)) = \frac{1}{(1-q)^{\chi(X)}}. 
\end{eqnarray*}
This is the proof of (\ref{gfEuler1}) in Zagier \cite{Zag}, \S 9.
Using (\ref{Euler1}) for $M^g$ with the action of $Z_g$,
Hirzebruch and H\"{o}fer \cite{Hir-Hof} showed that
\begin{eqnarray} \label{Euler21}
\chi(M, G) = \sum_{[g] \in G} \chi(M^g/Z_g).
\end{eqnarray}
Combining this expression with (\ref{fix2})
and using (\ref{gfEuler1}),
one gets
\begin{eqnarray*}
&& \sum_{n \geq 0} \chi(X^n, S_n) q^n \\
& = & \sum_{n \geq 0} q^n \sum_{\sum lN_l = n} \chi(\prod_{l=1}^n X^{(N_l)})
= \sum_{n \geq 0} 
\sum_{\sum lN_l = n} \prod_{l=1}^n (\chi(X^{(N_l)}) q^{lN_l}) \\
& = & \prod_{l \geq 1} \sum_{N_l \geq 0} \chi(X^{(N_l)}) q^{lN_l}
= \prod_{l \geq 1} \frac{1}{(1 - q^l)^{\chi(X)}}.
\end{eqnarray*}
This calculation was given by Hirzebruch and H\"{o}fer \cite{Hir-Hof} for 
surfaces.
It clearly works in general.

Now we come to the calculations for Poincar\'{e} polynomials.
Using the identification (\ref{fix2}),
one sees that
\begin{eqnarray*}
H^*(X^n, S_n) \cong \bigoplus_{\sum l N_l = n} H^*(\prod_{l=1}^n X^{(N_l)})
= \bigoplus_{\sum l N_l = n} \bigotimes_{N_l \geq 1} H^*(X^{(N_l)}).
\end{eqnarray*}
By the isomorphism (\ref{orbifoldcohomology}),
we have
\begin{eqnarray*}
H^*(X^n, S_n) \cong \bigoplus_{\sum l N_l = n} \bigotimes_{l=1}^n S^{N_l}(H^*(X)),
\end{eqnarray*}
as vector spaces.
This is what Vafa and Witten \cite{Vaf-Wit} used to prove (\ref{gfEuler2}).
We now analyze the grading shifts of the twisted sectors.
This can be reduced to the case of the action of an $n$-cycle $\sigma$ on $X^n$.
We assume that $X$ is a complex manifold for the time being.
By Example \ref{example:shift},
the shift is exactly the half of the real codimension of the diagonal.
This suggests that for any manifold $X$ of dimension $2m$,
not necessarily complex,
each $H^*(M^g)^{Z_g}$ should be shifted by 
half of the codimension of $(X^n)^g$ in $X^n$.
I.e., 
as $\bZ$-graded vector spaces,
\begin{eqnarray} \label{isomorphism}
H^*(X^n, S_n) \cong 
\bigoplus_{\sum l N_l = n} \bigotimes_{l=1}^n S^{N_l}(H^*(X)[m(l-1)]).
\end{eqnarray}
Then we have
\begin{eqnarray*}
&& \sum_{n \geq 0} P_t(X^n, S_n) q^n 
= \sum_{n \geq 0} q^n \sum_{\sum lN_l = n} 
P_t(\prod_{l = 1}^n S^{N_l}(H^*(X)[m(l-1)])) \\
& = & \sum_{n \geq 0} q^n \sum_{\sum lN_l = n} 
\prod_{l \geq 1} P_t(S^{N_l}(H^*(X)[m(l-1)])) \\
& = & \prod_{l \geq 1} \sum_{N_l \geq 0} P_t(S^{N_l}(H^*(X)[m(l-1)])) q^{lN_l} 
= \prod_{l \geq 1} \frac{\prod_{d \; odd} (1 + t^dq^l)^{b_{d - m(l-1)} (X)}}
{\prod_{d \; even} (1 - t^d q^l)^{b_{d - m(l - 1)}(X)}}.
\end{eqnarray*}
When $m$ is even,
we also have
\begin{eqnarray} \label{gfPoincare3}
&& \sum_{n \geq 0} P_t(X^n, S_n) q^n 
= \prod_{l \geq 1} \frac{\prod_{d \; odd} (1 + t^{m(l - 1) + d}q^l)^{b_d(X)}}
{\prod_{d \; even} (1 - t^{m(l - 1) + d}q^l)^{b_d(X)}}.
\end{eqnarray}

To summarize, 
we have proved the following 

\begin{theorem} \label{thm:main1}
For a $2m$-dimensional manifold $X$,
if $H^*(X^n, S_n)$ is graded as in (\ref{isomorphism}),
then the generating functional of the Poincar\'{e} polynomials of
$H^*(X^n, S_n)$ is given by:
\begin{eqnarray} \label{gfPoincare2}
\sum_{n \geq 0} P_t(X^n, S_n) q^n 
= \prod_{l \geq 1} \frac{\prod_{d \; odd} (1 + t^dq^l)^{b_{d - m(l-1)} (X)}}
{\prod_{d \; even} (1 - t^d q^l)^{b_{d - m(l - 1)}(X)}}.
\end{eqnarray}
\end{theorem}

\section{Hodge polynomials of the symmetric products}

Let $M$ be a complex manifold,
$G$ a finite group of bi-holomorphic transformations.
Then $M/G$ is a complex $V$-manifold,
hence one can consider the Dolbeault cohomology of $G$-invariant forms:
$$H^{*, *}(M/G) = H^{*, *}(M)^G.$$
Now we assume that $X$ is a complex manifold of complex dimension $m$.
Then the natural action of $S_n$ on $X^n$ is bi-holomorphic.
We then have
$$H^{*, *}(X^{(n)}) \cong S^n(H^{*, *}(X)).$$
By results from \S \ref{sec:series},
we then have the following

\begin{theorem} \label{thm:orbifoldHodge1}
For a compact complex manifold $X$,
we have
\begin{eqnarray} \label{gfHodge1}
\sum_{n \geq 0} h_{x, y}(X^{(n)}) q^n 
= \frac{\prod_{s+t \; odd} (1 + x^sy^t q)^{h^{s, t}(X)}}
{\prod_{s+t \; even} (1 - x^sy^t q)^{h^{s,t}(X)}}.
\end{eqnarray}
\end{theorem}

Recall Hirzebruch's $\chi_y$-genus for a compact complex manifold $M$ is
$$\chi_y(M) = \sum_{p, q \geq 0} (-1)^qh^{p, q}(M) y^p.$$
I.e. $\chi_y(M) = h_{y, -1}(M)$.
It is well-known that $\chi_{-1}(M) = \chi(M)$,
$\chi_1(M) = \sign(M)$ and $\chi_0(M) = p_a(M)$, the arithmetic genus of $M$
in Hirzebruch's sense.

\begin{corollary} \label{cor:orbifoldgenera1}
For a compact complex manifold $X$
we have
\begin{eqnarray} 
&& \sum_{n \geq 0} \chi_{-y}(X^{(n)}) q^n 
= \exp \left(\sum_{m \geq 1} \frac{\chi_{-y^m}(X)}{m} q^m\right), \label{gfchiy1} \\
&& \sum_{n \geq 0} p_a(X^{(n)}) q^n = \frac{1}{(1-q)^{p_a(X)}}, 
\label{gfarithmetic1}\\
&& \sum_{n \geq 0} \sign(X^{(n)}) q^n
= \frac{1}{(1-q^2)^{\chi(X)/2}}\left(\frac{1+q}{1-q}\right)^{\sign(X)/2}.
\label{gfsign1}
\end{eqnarray}
\end{corollary}

\begin{proof}
From (\ref{gfHodge1}),
it is easy to see that
\begin{eqnarray*} 
&& \sum_{n \geq 0} \chi_{-y}(X^{(n)}) q^n 
= \frac{\prod_{s+t \; odd} (1 + (-y)^s(-1)^t q)^{h_{s, t}(X)}}
{\prod_{s+t \; even} (1 - (-y)^s(-1)^t q)^{h_{s,t}(X)}} \\
& = & \frac{\prod_{s+t \; odd} (1 - y^s q)^{h_{s, t}(X)}}
{\prod_{s+t \; even} (1 - y^s q)^{h_{s,t}(X)}}
= \exp \left(-\sum_{s, t \geq 0} 
(-1)^{s+t} h^{s, t}(X) \log (1 - y^s q) \right) \\
& = & \exp \left(\sum_{s, t \geq 0} (-1)^{s+t} h^{s, t}(X)
\sum_{m>1} \frac{1}{m} (y^sq)^m \right) \\
& = & \exp \left(\sum_{m>1} \frac{1}{m} 
\sum_{s, t \geq 0} (-1)^{s+t} h^{s, t}(X) y^{ms}q^m \right) 
= \exp \left(\sum_{m>1} \frac{\chi_{-y^m}(X)}{m} q^m\right)
\end{eqnarray*}
(\ref{gfarithmetic1}) follows easily from (\ref{gfchiy}).
\begin{eqnarray*}
&& \sum_{n\geq 0} \sign(X^{(n)} q^n 
= \exp \left(\sum_{m>1} \frac{\chi_{-(-1)^m}(X)}{m} q^m\right) \\
& = & \exp \left(\sum_{l >1} \frac{\chi_{-(-1)^{2l}}(X)}{2l} q^{2l}\right) \cdot
\exp \left(\sum_{l >1} \frac{\chi_{-(-1)^{2l-1}}(X)}{2l-1} q^{2l-1}\right) \\
& = & \exp \left(\sum_{l >1} \frac{\chi(X)}{2l} q^{2l}\right) \cdot
\exp \left(\sum_{l >1} \frac{\sign(X)}{2l-1} q^{2l-1}\right) \\
& = & \frac{1}{(1-q^2)^{\chi(X)/2}}\left(\frac{1+q}{1-q}\right)^{\sign(X)/2}.
\end{eqnarray*}
\end{proof}

For a Riemann surface $X$,
the formulas in Theorem \ref{thm:orbifoldHodge1} 
and Corollary \ref{cor:orbifoldgenera1}
match with the results of Macdonald \cite{Mac2}.
Hirzebruch has computed the generating functional of the signature of
symmetric products of any compact oriented manifold
(for details, see Zagier \cite{Zag}, \S 9).
Formula (\ref{gfsign1}) matches with his result for complex manifolds.

Given a complex orbifold $M/G$, 
$M^g$ is a complex submanifold for every $g \in M$.
One can then consider the delocalized equivariant Dolbeault cohomology
$$H^{*, *}(M, G) = \bigoplus_{[g] \in G_*} H^{*, *}(M^g)^{Z_g}.$$
It follows from (\ref{fix2}) that
\begin{eqnarray*}
H^{*,*}(X^n, S_n) 
\cong \bigoplus_{\sum l N_l = n} H^{*,*}(\prod_{l=1}^n X^{(N_l)})
= \bigoplus_{\sum l N_l = n} \bigotimes_{l= 1}^n H^{*,*}(X^{(N_l)}).
\end{eqnarray*}
Using the isomorphism $H^{*, *}(M/G) \cong H^{*, *}(M)^G$,
we get
\begin{eqnarray*}
H^{*, *}(X^n, S_n) 
\cong \bigoplus_{\sum l N_l = n} \bigotimes_{l=1}^n S^{N_l}(H^{*, *}(X))
\end{eqnarray*}
as vector spaces.
Let $\dim_{\bC} X = 2k$ for some $k \in \frac{1}{2}\bZ$,
by Example \ref{example:shift}, 
as $\frac{1}{2}\bZ  \times \frac{1}{2}\bZ$-graded vector spaces,
\begin{eqnarray} \label{isomorphismD}
H^{*,*}(X^n, S_n) \cong 
\bigoplus_{\sum l N_l = n} \bigotimes_{l=1}^n 
	S^{N_l}(H^{*,*}(X)[k(l-1), k(l-1)]).
\end{eqnarray}
Of course, 
when $M$ is even dimensional,
both sided are $\bZ \times \bZ$-graded.
By (\ref{gfshifted}), 
we have
\begin{eqnarray*}
&& \sum_{n \geq 0} h_{x, y}(X^n, S_n) q^n \\
& = & \sum_{n \geq 0} q^n \sum_{\sum lN_l = n} 
h_{x, y}(\bigotimes_{l = 1}^n S^{N_l}(H^{*,*}(X)[k(l-1), k(l-1)])) \\
& = & \sum_{n \geq 0} q^n \sum_{\sum lN_l = n} 
\prod_{l \geq 1} h_{x, y}(S^{N_l}(H^{*,*}(X)[k(l-1), k(l-1)])) \\
& = & \prod_{l \geq 1} \sum_{N_l \geq 0} 
h_{x, y}(S^{N_l}(H^*(X)[k(l-1), k(l-1)])) q^{lN_l} \\
& = & \prod_{l \geq 1} 
\frac{\prod_{s+t \; odd} (1 + x^sy^tq^l)^{h^{s - k(l-1), t - k(l-1)} (X)}}
{\prod_{s+t \; even} (1 - x^sy^t q^l)^{h^{s - k(l - 1), t - k(l-1)}(X)}} \\
& = & \prod_{l \geq 1} 
\frac{\prod_{s+t \; odd} (1 + x^{s+k(l-1)}y^{t+k(l-1)}q^l)^{h^{s, t} (X)}}
{\prod_{s+t \; even} (1 - x^{s+k(l-1)}y^{t+k(l-1)} q^l)^{h^{s, t}(X)}}.
\end{eqnarray*}
To summarize, 
we have proved the following 

\begin{theorem} \label{thm:orbifoldHodge2}
For a closed complex manifold $X$ of dimension $2k$ 
for some $k \in \frac{1}{2}\bZ$,
if $H^{*,*}(X^n, S_n)$ is $\frac{1}{2}\bZ \times \frac{1}{2}\bZ$-graded 
as in (\ref{isomorphism}),
then we have
\begin{eqnarray} \label{gfHodge}
\sum_{n \geq 0} h_{x, y}(X^n, S_n) q^n 
= \prod_{l \geq 1} 
\frac{\prod_{s+t \; odd} (1 + x^{s+k(l-1)}y^{t+k(l-1)}q^l)^{h^{s, t} (X)}}
{\prod_{s+t \; even} (1 - x^{s+k(l-1)}y^{t+k(l-1)} q^l)^{h^{s, t}(X)}}.
\end{eqnarray}
\end{theorem}

For a complex orbifold $M/G$,
it is natural to define 
$$\chi_y(M, G) = \sum_{s, t \geq 0} (-1)^th^{s, t}(M, G)y^s.$$
Clearly we have $\chi_{-1}(M, G) = \chi(M, G)$. 
It is natural to regard $\chi_1(M, G)$ as the signature of $(M, G)$.
Similarly,
we regard $p_a(M, G) = \chi_0(M, G)$ as the arithmetic genus of $(M, G)$.
Similar to $\chi(M, G)$,
$\chi_y(M, G)$ is a sum of contributions from $M^g$, $[g] \in G_*$.
Denote by $h^{r, s}(M^g)^{Z_g}$ the dimension of $H^{r, s}(M^g)^{Z_g}
\cong H^{r, s}(M^g/Z_g)$.
Since
$$H^{p, q}(M, G) = \bigoplus_{[g] \in G_*} H^{p-F_g, q-F_g}(M^g)^{Z_g}
= \bigoplus_{[g] \in G_*} H^{p-F_g, q-F_g}(M^g/Z_g),$$
we have
\begin{eqnarray*}
\chi_y(M, G) & = & \sum_{[g] \in G_*} 
	\sum_{r, s \geq 0} (-1)^{t+F_g}h^{r, s}(M^g)^{Z_g} y^{s + F_g} \\
& = & \sum_{[g] \in G_*} 
	\sum_{r, s \geq 0} (-1)^{t+F_g}h^{r, s}(M^g/Z_g) y^{s + F_g}.
\end{eqnarray*}
Therefore,
\begin{eqnarray} \label{def:orbifoldchiy}
\chi_y(M, G) = \sum_{[g] \in G_*} (-y)^{F_g} \chi_y(M^g)^{Z_g}
= \sum_{[g] \in G_*} (-y)^{F_g} \chi_y(M^g/Z_g),
\end{eqnarray}
where
\begin{eqnarray*}
&& \chi_y(M^g)^{Z_g} = \sum_{s, t \geq 0} (-1)^th^{s, t}(M^g)^{Z_g}y^s, \\
&& \chi_y(M^g/Z_g) = \sum_{s, t \geq 0} (-1)^th^{s, t}(M^g/Z_g)y^s.
\end{eqnarray*}
Taking $y = 1$,
we get
\begin{eqnarray} \label{def:orbifoldsignature}
\chi_1(M, G) = \sum_{[g] \in G_*} (-1)^{F_g} \chi_y(M^g)^{Z_g}
= \sum_{[g] \in G_*} (-1)^{F_g} \chi_y(M^g/Z_g).
\end{eqnarray}
This indicates how one should define the orbifold signature.

\begin{corollary}
Under the assumptions of Theorem \ref{thm:orbifoldHodge2},
we have
\begin{eqnarray} 
&& \sum_{n \geq 0} \chi_{-y} (X^n, S_n) q^n 
=\exp \left( \sum_{n > 0} \frac{q^n}{n} \frac{\chi_{-y^n}(X)}{1- (y^kq)^n}
\right), \label{gfchiy} \\
&& \sum_{n\geq 0} p_a(X^n, S_n)q^n 
= \frac{1}{(1 - q)^{p_a(X)}}. \label{gfarithmetic2}
\end{eqnarray}
\end{corollary}

\begin{proof}
From (\ref{gfHodge}) we get
\begin{eqnarray*}
&& \sum_{n \geq 0} \chi_{-y}(X^n, S_n) q^n 
= \prod_{l \geq 1} 
\frac{\prod_{s+t \; odd} (1 + (-y)^{s+k(l-1)}(-1)^{t+k(l-1)}q^l)^{h^{s, t} (X)}}
{\prod_{s+t \; even} (1 - (-y)^{s+k(l-1)}(-1)^{t+k(l-1)} q^l)^{h^{s, t}(X)}} \\
& = & \prod_{l \geq 1}
\frac{\prod_{s+t \; odd} (1 - y^{s+k(l-1)} q^l)^{h^{s, t}(X)}}
{\prod_{s+t \; even} (1 - y^{s+k(l-1)} q^l)^{h^{s, t}(X)}}.
\end{eqnarray*}
The logarithm of the last term is
\begin{eqnarray*}
&& - \sum_{l \geq 1} \sum_{s, t \geq 0} 
(-1)^{s+t}h^{s, t}(X) \log(1 - y^{s+k(l-1)} q^l) \\
& = & \sum_{l \geq 1} \sum_{s, t \geq 0} (-1)^{s+t} h^{s, t}(X) 
\sum_{n > 0} \frac{1}{n} (y^{s+k(l-1)} q^l)^n \\
& = & \sum_{n > 0} \frac{1}{n} \sum_{l \geq 1} 
\sum_{s, t \geq 0} (-1)^{s+t} h^{s, t}(X) y^{sn+kn(l-1)} q^{ln}\\
& = &\sum_{n > 0} \frac{1}{n}\chi_{-y^n}(X) \sum_{l \geq 1} y^{kn(l-1)} q^{ln}
= \sum_{n > 0} \frac{q^n}{n} \frac{\chi_{-y^n}(X)}{1- (y^kq)^n}.
\end{eqnarray*}
When $x = 0$, $x^{s + k(l-1)} \neq 0$ only if $s = 0$ and $l = 1$,
hence
\begin{eqnarray*}
\sum_{n\geq 0} p_a(X^n, S_n)q^n 
= \frac{\prod_{t \; odd} (1+(-1)^tq)^{h^{0, t}(X)}}
{\prod_{t \; even} (1-(-1)^tq)^{h^{0, t}(X)}}
= \frac{1}{(1 - q)^{p_a(X)}}.
\end{eqnarray*}
\end{proof}

We omit the elementary proof of the following formulas:
\begin{eqnarray*}
&&\exp \left( \sum_{l \geq 1} 
\frac{q^{2l}}{2l} \frac{1}{1- q^{2l}} \right)
= \prod_{m \geq 1} \frac{1}{(1 - q^{2m})^{1/2}},  \\
&&  \exp \left( \sum_{l \geq 1} 
\frac{q^{2l-1}}{2l-1} \frac{1}{1- q^{2l-1}} \right) 
= \prod_{m \geq 1} \left(\frac{1 + q^m}{1-q^m}\right)^{1/2}, \\
&&  \exp \left( \sum_{l \geq 1} 
\frac{q^{2l-1}}{2l-1} \frac{1}{1+ q^{2l-1}} \right) 
= \prod_{m \geq 1} \left(\frac{1 + q^m}{1-q^m}\right)^{(-1)^{m+1}/2}.
\end{eqnarray*}

\begin{corollary}
Under the assumptions of Theorem \ref{thm:orbifoldHodge2},
we have
\begin{eqnarray*}
\sum_{n \geq 0} \chi_1(X^n, S_n) q^n 
=  \prod_{m \geq 1} \left(
\frac{1}{(1-q^{2m})^{\chi(X)/2}} 
\left(\frac{1+q^m}{1-q^m}\right)^{(-1)^{k(m+1)}\sign(X)/2}
\right).
\end{eqnarray*}
\end{corollary}

\begin{proof}
We take $y = -1$ in (\ref{gfchiy}).
For $k$ even, we have
\begin{eqnarray*}
&& \sum_{n \geq 0} \chi_1(X^n, S_n) q^n 
= \exp \left( \sum_{n > 0} 
\frac{q^n}{n} \frac{\chi_{-(-1)^n}(X)}{1- q^n} \right) \\
& = & \exp \left( \sum_{l \geq 1} 
\frac{q^{2l}}{2l} \frac{\chi_{-1}(X)}{1- q^{2l}} \right)\cdot
 \exp \left( \sum_{l \geq 1} 
\frac{q^{2l-1}}{2l-1} \frac{\chi_1(X)}{1- q^{2l-1}} \right) \\
& = & \prod_{m \geq 1} \left(
\frac{1}{(1-q^{2m})^{\chi(X)/2}} \left(\frac{1+q^m}{1-q^m}\right)^{\sign(X)/2}
\right).
\end{eqnarray*}
For $k$ odd, we have
\begin{eqnarray*}
&& \sum_{n \geq 0} \chi_1(X^n, S_n) q^n 
= \exp \left( \sum_{n > 0} 
\frac{(-1)^n}{n} \frac{\chi_{-(-1)^n}(X)}{1- (-q)^n} \right) \\
& = & \exp \left( \sum_{l \geq 1} 
\frac{q^{2l}}{2l} \frac{\chi_{-1}(X)}{1- q^{2l}} \right)\cdot
 \exp \left( \sum_{l \geq 1} 
\frac{q^{2l-1}}{2l-1} \frac{\chi_1(X)}{1 + q^{2l-1}} \right) \\
& = & \prod_{m \geq 1} \left(
\frac{1}{(1-q^{2m})^{\chi(X)/2}} 
\left(\frac{1+q^m}{1-q^m}\right)^{(-1)^{m+1}\sign(X)/2}
\right).
\end{eqnarray*}
\end{proof}

\section{The $B$-version}

In the study of mirror symmetry of Calabi-Yau manifolds,
one encounters another interesting Dolbeault cohomology algebra
of a closed complex manifold $M$:
$$H^{-*, *}(M) = \bigoplus_{p, q \geq 0} H^q(M, \Lambda^p(TM)) 
= \bigoplus_{p, q \geq 0} H^{-p, q}(M),$$
with the algebra structure induced from the wedge products.
We will refer to it as the {\em $B$-algebra}.
Clearly all the results in last section have a version for $H^{-*, *}$,
so we will be brief.

Set $h^{-p, q}(M) = \dim H^{-p, q}(M)$.
The $B$-Hodge polynomial is by definition:
$$\hat{h}_{x, y}(M) = \sum_{p, q \geq 0} h^{-p, q}(M) x^py^q.$$
For a Calabi-Yau $d$-fold $M$,
Serre duality shows that $h^{-p, q}(M) = h^{d-p, q}$.
Hence we have
\begin{eqnarray*}
&& \hat{\chi}_{-y}(M) = \sum_{p, q \geq 0} (-1)^qh^{d-p, q}(M)(-y)^p \\
& = & \sum_{p, q \geq 0} (-1)^qh^{p, q}(M)(-y)^{d - p}
= (-y)^d \chi_{-y^{-1}}(M).
\end{eqnarray*}
Or equivalently
\begin{eqnarray}
y^{d/2}\hat{\chi}_{-y}(M) = (-1)^{d/2} \cdot y^{-d/2}\chi_{-y^{-1}}(M).
\end{eqnarray}
For a complex orbifold $M/G$,
$$H^{-*, *}(M/G) = H^{-*, *}(M)^G.$$
One define $\hat{h}_{x, y}(M/G)$ in a similar fashion.
In particular,
$$H^{*, *}(X^{(n)}) \cong S^n(H^{*, *}(X)),$$
for any closed complex manifold $X$.
One then has the following analogue of Theorem \ref{thm:orbifoldHodge1}:

\begin{theorem} \label{thm:orbifoldHodge1B}
For a compact complex manifold $X$,
we have
\begin{eqnarray} \label{gfHodge1B}
\sum_{n \geq 0} \hat{h}_{x, y}(X^{(n)}) q^n 
= \frac{\prod_{s+t \; odd} (1 + x^sy^t q)^{h^{-s, t}(X)}}
{\prod_{s+t \; even} (1 - x^sy^t q)^{h^{-s,t}(X)}}.
\end{eqnarray}
\end{theorem}

The $B$-version of Hirzebruch's $\chi_y$-genus 
for a compact complex manifold $M$ is defined to be
$$\hat{\chi}_y(M) = \sum_{p, q \geq 0} (-1)^qh^{-p, q}(M) y^p 
= \hat{h}_{y, -1}(M).$$

\begin{corollary} \label{cor:orbifoldgenera1B}
For a compact complex manifold $X$
we have
\begin{eqnarray} 
\sum_{n \geq 0} \hat{\chi}_{-y}(X^{(n)}) q^n 
= \exp \left(\sum_{m \geq 1} \frac{\hat{\chi}_{-y^m}(X)}{m} q^m\right).
	\label{gfchiy1B} 
\end{eqnarray}
\end{corollary}

Given a complex orbifold $M/G$, 
consider the $B$-delocalized equivariant Dolbeault cohomology
$$H^{-*, *}(M, G) = \bigoplus_{[g] \in G_*} H^{-*, *}(M^g)^{Z_g}.$$
We give it a bigrading by fractional numbers by
$$H^{p, q}(M, G) = \bigoplus_{[g] \in G_*} H^{p-F_g, q-F_g}(M^g)^{Z_g}
= \bigoplus_{[g] \in G_*} H^{p-F_g, q-F_g}(M^g/Z_g).$$
In particular,
it follows from (\ref{fix2}) and Example \ref{example:shift} that
\begin{eqnarray} \label{isomorphismB}
H^{-*,*}(X^n, S_n) \cong 
\bigoplus_{\sum l N_l = n} \bigotimes_{l=1}^n S^{N_l}(H^{-*, *}(X)[k(l-1), k(l-1)]).
\end{eqnarray}
where $\dim_{\bC} X = 2k$ for some $k \in \frac{1}{2}\bZ$.
Again,
when $M$ is even dimensional,
both sided are $\bZ \times \bZ$-graded.

\begin{theorem} \label{thm:orbifoldHodge2B}
For a closed complex manifold $X$ of dimension $2k$ 
for some $k \in \frac{1}{2}\bZ$,
if $H^{*,*}(X^n, S_n)$ is $\frac{1}{2}\bZ \times \frac{1}{2}\bZ$-graded 
as in (\ref{isomorphismB}),
then we have
\begin{eqnarray} \label{gfHodgeB} \;\;\;\;\;\;\;
\sum_{n \geq 0} \hat{h}_{x, y}(X^n, S_n) q^n 
= \prod_{l \geq 1} 
\frac{\prod_{s+t \; odd} (1 + x^{s+k(l-1)}y^{t+k(l-1)}q^l)^{h^{-s, t} (X)}}
{\prod_{s+t \; even} (1 - x^{s+k(l-1)}y^{t+k(l-1)} q^l)^{h^{-s, t}(X)}}.
\end{eqnarray}
\end{theorem}

Define 
$$\hat{\chi}_y(M, G) = \sum_{s, t \geq 0} (-1)^th^{-s, t}(M, G)y^s.$$
Then we have

\begin{corollary}
Under the assumptions of Theorem \ref{thm:orbifoldHodge2B},
we have
\begin{eqnarray} 
\sum_{n \geq 0} \hat{\chi}_{-y} (X^n, S_n) q^n 
=\exp \left( \sum_{n > 0} \frac{q^n}{n} \frac{\hat{\chi}_{-y^n}(X)}{1- (y^kq)^n}
\right). \label{gfchiyB} 
\end{eqnarray}
\end{corollary}

Similar to $\chi(M, G)$ and $\chi_y(M, G)$,
we have
\begin{eqnarray} \label{def:orbifoldchiyB}
\hat{\chi}_y(M, G) 
= \sum_{[g] \in G_*} (-y)^{F_g} \hat{\chi}_y(M^g/Z_g).
\end{eqnarray}

\section{A special case of a formula for elliptic genera of symmetric products}

Our results above give a mathematical proof of a special case of a formula for
the elliptic genera of symmetric products found
by Dijkgraaf, Moore, E. Verlinde and H. Verlinde \cite{Dij-Moo-Ver-Ver}.
For a closed complex $d$-manifold $M$, 
expand
\begin{eqnarray*}
E_{q, y} = y^{-\frac{d}{2}} \bigotimes_{n \geq 1} 
(\Lambda_{-yq^{n-1}}TM \otimes \Lambda_{-y^{-1}q^n}T^*M \otimes
S_{q^n}TM \otimes S_{q^n}T^*M)
\end{eqnarray*}
as a sum of holomorphic vector bundles:
$$E_{q, y} = \bigoplus_{m \geq 0, l} q^my^l E_{m. l}.$$
It is easy to see that each $E_{l, m}$ is of finite rank. 
Let $c(m, l)$ be the Riemann-Roch number of $E_{m, l}$, 
then the elliptic genus defined in \cite{Dij-Moo-Ver-Ver} is
\begin{eqnarray*}
\chi(M; q, y) = \sum_{m \geq 0, l} c(m, l)q^m y^l.
\end{eqnarray*}
It was also described in \cite{Dij-Moo-Ver-Ver} as the trace of some operator
on some Hilbert space $\cH(M)$,
but the exact description of this Hilbert space and the operators involved
were not explicitly given. 
Furthermore,
for a complex orbifold $M/G$,
the definition of its orbifold elliptic genus used such a description:
the total orbifold Hilbert space takes the form 
$$\cH(M, G) = \oplus_{[g] \in G_*} \cH_g^{Z_g},$$
where $\cH_g$ is a Hilbert space with a $Z_g$ action on it.
String theoretical arguments used in \cite{Dij-Moo-Ver-Ver}
requires the $Z_g$-action satisfies some nice properties in the case of
symmetric products,
and the following formula was derived
\begin{eqnarray} \label{gfelliptic}
\sum_{n \geq 0} p^n \chi(X^n, S_n;q, y)
= \prod_{n>0, m \geq 0, l} \frac{1}{(1-p^nq^my^l)^{c(nm, l)}},
\end{eqnarray}
where $\chi(X^n, S_n;q, y)$ denotes the orbifold elliptic genus of the 
symmetric product, $\chi(X; q, y) = \sum c(m, l) q^m y^l$.

Mathematically,
the definition of the orbifold elliptic genus has yet to be worked out.
But the following discussion provides some hint on the right definition.
Presumably,
it involves the contributions from the fixed point sets $M^g$ for $[g]$
running through all conjugacy classes.
We look at two special cases.
When $y = 1$,
Dijkgraaf {\em et al} \cite{Dij-Moo-Ver-Ver} claimed that
$\chi(M; q, 1) = \chi(M)$, the Euler number of $M$.
The orbifold Euler number as shown by Hirzebruch and H\"{o}fer \cite{Hir-Hof}
is given by formula (\ref{orbifoldEuler2}) in \S 1:
$$\chi(M, G) = \sum_{[g] \in G_*} \chi(M^g/Z_g).$$
When $q = 0$,
Dijkgraaf {\em et al} \cite{Dij-Moo-Ver-Ver} claimed that
\begin{eqnarray*}
\chi(M;0, y) = \sum_{r, s} (-1)^{r+s}y^{r-\frac{d}{2}}h^{r, s}(M) 
= y^{-\frac{d}{2}}\chi_{-y}(M).
\end{eqnarray*}
(It seems to the author that $\chi(M;0, y)$ should be equal to 
$y^{-\frac{d}{2}}\hat{\chi}_{-y}(M)$,
since $E_{0, y} = y^{-\frac{d}{2}}\Lambda_{-y}TM$.)
Formula (\ref{gfelliptic}) reduces to
\begin{eqnarray} \label{gfq=01}
\sum_{n \geq 0} p^n \chi(X^n, S_n;0, y)
= \prod_{n>0, l\geq 0} \frac{1}{(1-p^ny^l)^{c(0, l)}}.
\end{eqnarray}
We now rewrite the right hand side as follows:
\begin{eqnarray*}
&& \prod_{n>0, l\geq 0} \frac{1}{(1-p^ny^l)^{c(0, l)}}
= \exp \left( \sum_{n > 0} \sum_{l \geq 0} -c(0, l) \log (1 - p^ny^l) \right)\\
& = & \exp \left(\sum_{n > 0} \sum_{l \geq 0} c(0, l) 
	\sum_{m> 0} \frac{(p^ny^l)^m}{m}\right)
= \exp \left(\sum_{n > 0} \sum_{m> 0} \chi(X; 0, y^m)\frac{p^{nm}}{m} \right)\\
& = & \exp \left(\sum_{m> 0} \frac{\chi(X; 0, y^m)}{m} \frac{p^m}{1 - p^m}\right).
\end{eqnarray*}
So (\ref{gfq=01}) is equivalent to 
\begin{eqnarray} \label{gfq=02}
\sum_{n \geq 0} p^n \chi(X^n, S_n;0, y)
= \exp \left(\sum_{m> 0} \frac{\chi(X; 0, y^m)}{m} \frac{p^m}{1 - p^m}\right).
\end{eqnarray}
Noting the similarity between (\ref{gfq=02}) and (\ref{gfchiy}),
we are led to the following:

\begin{theorem}
If one defines 
\begin{eqnarray} \label{def:q=0}
\chi(M, G;0, y) = y^{-\frac{\dim M}{2}}\chi_{-y}(M, G)
\end{eqnarray}
for complex orbifolds $M/G$,
then (\ref{gfq=02}) holds for any closed complex manifold $X$.
\end{theorem}

\begin{proof}
We use (\ref{gfchiy}) for $q= y^{-k}p$,
where $\dim_{\bC} = 2 k$ for some $k \in \frac{1}{2}\bZ$:
\begin{eqnarray*}
&& \sum_{n \geq 0} p^n \chi(X^n, S_n;0, y)
=\sum_{n \geq 0} \chi_y(X^n, S_n) y^{-kn}p^n \\
& = & \exp \left( \sum_{m > 0} \frac{(y^{-k}p)^m}{m} 
	\frac{\chi_{-y^m}(X)}{1-p^m}  \right) \\
& = & \exp \left( \sum_{m > 0} \frac{(y^m)^{-k}\chi_{-y^m}(X)}{m} 
	\frac{p^m}{1-p^m}  \right) \\
& = & \exp \left(\sum_{m> 0} 
	\frac{\chi(X; 0, y^m)}{m} \frac{p^m}{1 - p^m}\right).
\end{eqnarray*}
\end{proof}

By (\ref{def:orbifoldchiy}),
(\ref{def:q=0}) becomes
\begin{eqnarray*}
\chi(M, G; 0, y) 
& = & y^{-\frac{\dim M}{2}} \sum_{[g] \in G_*} y^{F_g} \chi_{-y}(M^g/Z_g) \\
& = & \sum_{[g] \in G_*} y^{-\frac{\dim M}{2} + F_g + \frac{\dim M^g}{2}} 
  \chi(M^g/Z_g; 0, y).
\end{eqnarray*}
In the case of symmetric products, 
$F_g = \frac{\dim M - \dim M^g}{2}$,
hence if we use (\ref{def:q=0}), then
\begin{eqnarray*}
\chi(X^n, S_n; 0, y) = \sum_{[g] \in (S_n)_*} \chi((X^n)^g/Z_g; 0, y).
\end{eqnarray*}
It seems reasonable to conjecture that
if one defines
$$\chi(X^n, S_n; q, y) = \sum_{[g] \in (S_n)_*} \chi((X^n)^g/Z_g; q, y),$$
then (\ref{gfelliptic}) holds.

We also have the $B$-version of the above results:

\begin{theorem}
If one defines 
\begin{eqnarray*} 
\hat{\chi}(M, G;0, y) = y^{-\frac{\dim M}{2}}\hat{\chi}_{-y}(M, G)
\end{eqnarray*}
for complex orbifolds $M/G$,
then 
\begin{eqnarray*}
\hat{\chi}(X^n, S_n; 0, y) 
= \sum_{[g] \in (S_n)_*} \hat{\chi}((X^n)^g/Z_g; 0, y).
\end{eqnarray*}
Furthermore, (\ref{gfq=02}) holds with $\chi$ replaced by $\hat{\chi}$
for any closed complex manifold $X$.
\end{theorem}

\section{Heisenberg superalgebra representation}

Set $\cF(X) = \sum_{n \geq 0} H^*(X^n, S_n)$.
Vafa and Witten \cite{Vaf-Wit} noticed that this space admits an action
by an infinite dimensional Heisenberg superalgebra.
Following Segal \cite{Seg} and Wang \cite{Wan},
we give a construction of such an action 
in terms of the some functors defined in \S \ref{sec:delocalized}.
We first recall some well-known facts about Heisenberg superalgebras
(see e.g. Kac \cite{Kac}).
Let $\fh$ be a finite dimensional complex superspace 
with a non-degenerate supersymmetric 
bilinear form $\eta(\cdot, \cdot)$.
In other words,
$\fh$ is $\bZ_2$-graded: $\fh = \fh_{\bar{0}} \oplus \fh_{\bar{1}}$,
$\eta: \fh \otimes \fh \to \bC$ is a nondegenerate bilinar functional
of degree $\bar{0}$,
such that 
$$\eta(a, b) = (-1)^{\bar{a} \cdot \bar{b}} \eta(b, a),$$
for homogeneous $a, b \in \fh$.
Here we have used $\bar{a}$ to denote the degree of $a$.
We will refer to $\eta$ simply as a pairing on $\fh$.
Consider the Lie superalgebra (Kac \cite{Kac}, \S 3.5)
$$\hat{\fh} = \bC[t, t^{-1}] \otimes_{\bC} \fh \oplus \bC K$$
with commutation relations 
$$[a_m, b_n] = m\cdot \eta(a, b) \delta_{m, -n}K, \;\;\;\; [K, \hat{\fh}] = 0,$$
where $m, n \in \bZ, a, b \in \fh$, $a_m = t^m a$ and $b_n = t^n b$.
Clearly $\hat{\fh}$ is the direct sum of $\fh$ with Heisenberg
superalgebra $\hat{\fh}'$ given by
$$\hat{\fh}' = \hat{\fh}^{<} \oplus \bC K \oplus \hat{\fh}^{>},$$
where $\hat{h}^< = \oplus_{n \geq 1} (t^{-n}\fh)$ and
$\hat{\fh}^> = \oplus_{n \geq 1} (t^n\fh)$.
Let $\hat{\fh}^+ = \hat{\fh}^> \oplus \bC K$. 
Given $k \in \bC$,
denote by $\pi^k$ the $1$-dimensional representation of $\hat{\fh}^+$ defined by
$$\pi^k(\hat{\fh}) = 0, \;\;\;\; \pi(K) = k.$$
The Verma module $\widetilde{V}^k(\fh, \eta) = U(\hat{\fh}) \otimes_{U(\hat{\fh}^+)} \pi^k$
has an explicit description as follows:
\begin{eqnarray} \label{Verma}
\widetilde{V}^k(\fh, \eta) = S^*(\hat{\fh}^<) 
= S^*(\oplus_{l \geq 1} t^{-l} \fh) = \otimes_{l \geq 1} S^*(t^{-l}\fh),
\end{eqnarray}
where $K$ acts as multiplication by a constant $k$,
for $m > 0$, 
$a_{-m}$ acts by multiplication,
$a_m$ by $k m$ times the contraction:
$$(t^m \otimes a) \cdot (t^{-n} \otimes b) = k m \delta_{m, n} \eta(a, b).$$
When $k \neq 0$,
there is a structure of a simple conformal vertex algebra on $\widetilde{V}^k$
(Kac \cite{Kac}, \S 4.7 and \S 4.10).
It is called the {\em free bosonic vertex algebra}.
This vertex operator algebra structure is essentially determined
by the action of $\hat{\fh}'$ on $\widetilde{V}^k(\fh, \eta)$.

It is straightforward to generalize to the $\bZ$ or $\bZ \times \bZ$-graded
version.
For a $\bZ$-graded $\fh$,
we denote by $|a|$ the degree of a homogeneous element $a \in \fh$.
We assume that the $\eta$ is a nondegenerate graded symmetric bilinear form
of degree $0$,
we can give $t$ an arbitrary even degree $-d$. 
Then $K$ acts as an operator of degree $0$,
and for $m > 0$,
$a_{-m}$ acts by an operator of degree $|a| + md$,
$a_m$ acts by an operator of degree $- |a| - md$.
Similarly,
for a $\bZ \times \bZ$-graded $\fh$,
we  denote by $||a||$ the bi-degree of a homogeneous element $a \in \fh$.
We assume that $\eta$ is a nondegenerate bilinear form of bi-degree $(0, 0)$,
bigraded symmetric in the sense that
$$\eta(a, b) = (-1)^{(p_1 + q_1)(p_2, q_2)} \eta(b, a)$$
for $a \in \fh_{p_1. q_1}$, $b \in \fh_{p_2, q_2}$.
We give $t$ a bi-degree $(p, p)$.
Then $K$ acts as an operator of bi-degree $(0, 0)$,
and for $m > 0$, $a \in \fh_{r, s}$
$a_{-m}$ acts as an operator of bidegree $(r + mp, t + mp)$,
$a_m$ acts as an operator of bidegree $(-(r + mp), -(t + mp))$.

From the isomorphism (\ref{isomorphism}),
we get
\begin{eqnarray*}
\cF(X) & \cong & \sum_{n \geq 0} \bigoplus_{\sum l N_l = n} \bigotimes_{l=1}^n 
S^{N_l}(H^*(X)[d(l-1)]) \\
& = & \bigotimes_{l \geq 1} \sum_{N_l \geq 0} S^{N_l}(H^*(X)[d(l-1)])
= \bigotimes_{l \geq 1} S^*(H^*(X)[d(l-1)]) \\
& = & S^*(\oplus_{l\geq 1} H^*(X)[d(l-1)]).
\end{eqnarray*}
Comparing with (\ref{Verma}),
we see that when $d$ is an even number,
if $t$ has degree $-d$,
then we have
$$ \sum_{n \geq 0} H^*(X^n, S_n) \cong \widetilde{V}^k(H^*(X)[-d], \eta),$$
where $\eta(\alpha, \beta) = \int_X \alpha \cup \beta$ for 
$\alpha, \beta \in H^*(X)$.
Here we have assumed that $X$ is a connected closed oriented manifold 
of dimension $2d$ with $d$ even,
hence by Poincar\'{e} duality,
$\eta$ is a nondegenerate graded symmetric bilinear form on $H^*(X)[-d]$.

It follows from the above discussion that $\cF(X)$
is the Verma module of the Heisenberg superalgebra for 
$(H^*(X)[-d], \eta)$.
Now we give geometric realization of the relevant operators.
We first recall some maps defined in \S \ref{sec:delocalized}.
For any positive integers $p < n$,
the inclusion $S_p \times S_{n-p} \hookrightarrow S_n$ induces two 
maps
\begin{eqnarray*}
&& \Res: H^*(X^n, S_n) \to H^*(X^n, S_p \times S_{n-p}), \\
&& \Ind: H^*(X^n, S_p \times S_{n-p}) \to H^*(X^n, S_n).
\end{eqnarray*}
We also have the K\"{u}nneth isomorphism
$$\kappa: H^*(X^n, S_p \times S_{n-p}) \to
	H^*(X^p, S_p) \otimes H^*(X^{n-p}, S_{n-p}).$$
Clearly $\Res, \Ind$ and $\kappa$ all have degree $0$
with the given $\bZ$-grading on relevant spaces.
For $0 \leq m \leq n$,
define 
$$\cdot: H^r(X^m, S_m) \otimes H^*(X^{n-m}, S_{n-m}) \to H^*(X^n, S_n)$$
as the composition
\begin{eqnarray*}
H^*(X^m, S_m) \otimes H^*(X^{n-m}, S_{n-m}) 
\stackrel{\kappa^{-1}}{\longrightarrow} H^*(X^n, S_m \times S_{n-m})
\stackrel{\Ind}{\longrightarrow} H^*(X^n, S_n),
\end{eqnarray*}
also define
$$\Delta: H^*(X^n, S_n) \to \oplus_{m = 0}^n 
H^*(X^m, S_m) \otimes H^*(X^{n-m}, S_{n-m})$$
as the direct sum of the compositions
\begin{eqnarray*}
H^*(X^n, S_n) \stackrel{\Res}{\longrightarrow} H^*(X^n, S_m \times S_{n-m})
\stackrel{\kappa}{\longrightarrow} 
	H^*(X^m, S_m) \otimes H^*(X^{n-m}, S_{n-m}).
\end{eqnarray*}
Hence we get maps $\cdot: \cF(X) \otimes \cF(X) \to \cF(X)$ and 
$\Delta: \cF(X) \to \cF(X) \otimes \cF(X)$.
It is easy to verify that under the identification,
$\cF(X)  \cong S^*(\oplus_{l \geq 1} t^l H^*(X)[-m])$,
$\cdot$ and $\Delta$ gives the multiplication and comultiplication for
the standard Hopf algebra structure on the latter space.
For an $n$-cycles $\sigma_n \in S_n$, $Z_{\sigma_n} = \langle \sigma_n \rangle$
and $H^*((X^n)^{\sigma_n})^{Z_{\sigma_n}} \cong H^*(X)$.
So there are two homomorphisms
\begin{eqnarray*}
&& i_n: H^*(X) \to H^*(X^n, S_n), \\
&& \pi_n: H^*(X^n, S_n) \to H^*(X)
\end{eqnarray*}
induced by the inclusion and projection respectively.
When regarded as a map from $H^*(X)[-d]$ to $H^*(X^n, S_n)$,
$i_n$ has degree $nd$;
similarly,
$\pi_n$ is a map of degree $-nd$ from $H^*(X^n, S_n)$ to $H^*(X)[-d]$.
For any $1 \leq m \leq n$ and $\alpha \in H^*(X)$,
define 
$$p_m(\alpha): H^*(X^n, S_n) \to H^*(X^{n-m})$$
as $m$ times the composition 
\begin{eqnarray*}
H^*(X^n, S_n) & \stackrel{\Delta}{\longrightarrow} & 
	H^*(X^m, S_m) \otimes H^*(X^{n-m}, S_{n-m}) \\ 
& \stackrel{\pi_m \otimes 1}{\longrightarrow} &
	H^*(X) \otimes H^*(X^{n-m}, S_{n-m})
\stackrel{\iota_{\alpha} \otimes 1}{\longrightarrow} H^*(X^{n-m}, S_{n-m}).
\end{eqnarray*}
Also define 
$$q_m(\alpha): H^*(X^{n-m}) \to H^*(X^n, S_n)$$
as the composition 
\begin{eqnarray*}
H^*(X^{n-m}, S_{n-m}) 
& \stackrel{\alpha \otimes 1}{\longrightarrow} &
 	H^*(X) \otimes H^*(X^{n-m}, S_{n-m}) \\
& \stackrel{i_m \otimes 1}{\longrightarrow} &
 	H^*(X^m, S_m) \otimes H^*(X^{n-m}, S_{n-m}) 
\stackrel{\cdot}{\longrightarrow} H^*(X^n, S_n)
\end{eqnarray*}
Then $p_m(\alpha)$ has degree $- (|\alpha| - d) - md$,
and $q_m(\alpha)$ has degree $(|\alpha| - d) + md$. 
It is straightforward to check that
$p_m(\alpha)$ is $m$ times the contraction by $t^m \otimes \alpha$ 
and $q_m(\alpha)$ is the multiplication by $t^{-m} \otimes \alpha$.
For an open manifold $X$,
one can use the natural pairing on the direct sum of $H^*_c(X)$
with its dual space to construct a Heisenberg superalgebra and its action
on $\cF_c(X) = \sum_{n\geq 0} H^*_c(X^n, S_n)$.

When $M$ is a closed complex manifold,
one can also consider the space
$$\widehat{\cF}(X) = \sum_{n \geq 0} H^{-*, *}(X^n, S_n).$$
When $X$ is Calabi-Yau $d$-fold,
there is a pairing $\hat{\eta}$ given by the composition:
$$H^{-p, q}(X) \otimes H^{-r, s}(X) 
\stackrel{\wedge}{\rightarrow} H^{-(p+r), q+s}(X) 
\stackrel{\sharp}{\to} H^{d-(p+r), q+s}(X) 
\stackrel{\int_X}{\to} \bC,$$
where $\sharp$ is the isomorphism induced by the holomorphic volume form on $X$.
Similar to the above discussion,
one gets an action of the Heisenberg superalgebra constructed from
$(H^{-*, *}(X), \eta)$ on $\widehat{\cF}(X)$.
In general,
one use the natural pairing on $H^{-*, *}(M)$ and its dual space
to construct a Heisenberg superalgebra and its representation on
$\widehat{\cF}(X)$.

\section{Surface case}

When $X$ is a smooth algebraic surface,
the Hilbert scheme $X^{[n]}$ of $0$-dimensional subscheme of length $n$
is a smooth algebraic variety.
G\"{o}ttsche \cite{Got} has shown that
$$\sum_{n \geq 0} P_t(X^{[n]}) q^n 
= \prod_{l \geq 1} 
\frac{(1+t^{2l-1}q^l)^{b_1(X)}(1 + t^{2l+1}q^l)^{b_3(X)}}
{(1 - t^{2l-2}q^l)^{b_0(X)}(1 - t^{2l}q^l)^{b_2(X)}
(1 - t^{2l+2}q^l)^{b_4(X)}}.$$
Also G\"{o}ttsche-Soergel \cite{Got-Soe} and 
Cheah \cite{Che} have shown by different methods that
\begin{eqnarray*} 
\sum_{n \geq 0} h_{x, y}(X^{[n]}) q^n 
= \prod_{l \geq 1} 
\frac{\prod_{s+t \; odd} (1 + x^{s+(l-1)}y^{t+(l-1)}q^l)^{h^{s, t} (X)}}
{\prod_{s+t \; even} (1 - x^{s+(l-1)}y^{t+(l-1)} q^l)^{h^{s, t}(X)}}.
\end{eqnarray*}
Compared with our formulas for $H^*(X^n, S_n)$ and $H^{*, *}(X^n, S_n)$,
we see that for a smooth algebraic surface $X$,
\begin{eqnarray}
&& H^*(X^{[n]}) \cong H^*(X^n, S_n),  \label{HS1} \\
&& H^{*, *}(X^{[n]}) \cong H^{*, *}(X^n, S_n) \label{HS2}
\end{eqnarray}
as (bi-)graded vector spaces.
As a consequence, 
$\chi(X^{(n)}) = \chi(X^n, S_n)$ for smooth algebraic surfaces.
This was noticed by Hirzebruch and H\"{o}fer \cite{Hir-Hof}.
In \cite{Got-Soe},
the grading shift actually has an explanation in terms of intersection
cohomology.

The isomorphisms (\ref{HS1}) and (\ref{HS2}) are special cases 
of a very general result.
As was proved in a paper of Batyrev and Dais \cite{Bat-Dai}, 
$h^{p,q}(M,G) = h^{p,q}(\widehat{M/G})$ for 
a crepant, full desingularization $\widehat{M/G}\to M/G$ of the orbifold $M/G$, 
whenever the so-called ``strong McKay correspondence" 
holds for all quotient singularities occurring along every stratum of $X/G$. 
For a smooth algebraic surface,
the Hilbert scheme $X^{[n]}$ provides a crepant resolution of $X^{(n)}$.
G\"{o}ttsche \cite{Got2} showed that this resolution satisfies the 
``strong McKay correspondence",
while Batyrev proved it in general in his recent preprint \cite{Bat}.
If one uses such results,
our calculations give a method of computing the Hodge numbers of $X^{[n]}$.

It seems to be a general phenomenon that invariants of $X^{[n]}$ 
can be computed by suitably defined corresponding invariants of $X^{(n)}$ for
surfaces. 
If one can relate the construction of Nakajima \cite{Nak}
and Grojnowski \cite{Gro} for the Hilbert schemes to our construction in \S 7
or the constructions in Segal \cite{Seg} and Wang \cite{Wan},
then one obtains another example of this phenomenon.
Also related is a conjecture in 
Dijkgraaf-Moore-Verlinde-Verlinde \cite{Dij-Moo-Ver-Ver}
on the elliptic genera for the Hilbert schemes of K3 surfaces. 
Liu \cite{Liu} has proposed a method to directly verify it.
Nevertheless the study of the relationship with 
the elliptic genera of the symmetric
products still seems interesting.

{\bf Acknowledgement}.
{\em The work in this paper is carried out during the author's visit at Texas
A$\&$M University. 
He thanks the Mathematics Department and the GAT group for hospitality 
and financial support.
The author likes to thank Huai-Dong Cao, Weiping Li, 
Kefeng Liu and Weiqiang Wang
for helpful discussions and encouragement.
Finally some email exchanges with Cumrun Vafa helped 
to understand the arguments for symmetric products in \cite{Vaf-Wit}.
}

\end{document}